\documentclass[]{amsart}
\usepackage{amssymb}
\usepackage[all]{xy}% xymatrix for CD's
\newdir{ >}{{}*!/-10pt/@{>}}% tail for mono arrow
\newdir^{ (}{{}*!/-5pt/@^{(}}% tail for hookarrow
\newdir_{ (}{{}*!/-5pt/@_{(}}% tail for hookarrow
\allowdisplaybreaks

\def\undersetbrace#1\to#2{\underbrace{#2}_{#1}}
\def\oversetbrace#1\to#2{\overbrace{#2}^{#1}}
\def\AMSunderset#1\to#2{\underset{#1}{#2}}
\def\AMSoverset#1\to#2{\overset{#1}{#2}}
\def\therosteritem#1{{\rm (#1)}}
\def\East#1#2{\overset{#1}{\longrightarrow}}
\swapnumbers

\newtheorem{prop}[subsection]{Proposition}
\newtheorem*{prop*}{Proposition}
\newtheorem{thm}[subsection]{Theorem}
\newtheorem*{thm*}{Theorem}
\newtheorem{lem}[subsection]{Lemma}
\newtheorem*{lem*}{Lemma}
\newtheorem{cor}[subsection]{Corollary}
\newtheorem*{cor*}{Corollary}

\newtheorem*{conj*}{Conjecture}
\newenvironment{demo}[1]{\par\smallskip\noindent{\bf #1.}}{\par\smallskip}

\def\fr{\frak}
\def\a{\alpha}
\def\al{\alpha}

\def\ga{\gamma}
\def\Ga{\Gamma}

\def\ra{\rightarrow}

\def\lra{\longrightarrow}

\def\lRa{\Longrightarrow}

\def\lBia{\Longleftrightarrow}

\def\ua{\uparrow}
\def\da{\downarrow}

\def\hra{\hookrightarrow}

\def\loopra{\looparrowright}

\def\bra{\mapsto}

\def\bsl{\backslash}

\def\ocong#1{\AMSoverset {#1} \to \cong}
\def\oset#1{\AMSoverset {#1} \to \ra}
\def\osetl#1{\AMSoverset {#1} \to \lra}

\def\set#1#2{ \{  {#1} \mid {#2} \} }

\def\abs#1{|{#1}|}

\def\Cinfty#1{C^{\infty}(#1)}

\def\ti{\tilde}
\def\wti{\widetilde}

\def\what{\widehat}
\def\wbar{\overline}
\def\vphi{\varphi}

\def\F{{\mathcal{F}}}
\def\G{{\mathcal{G}}}
\def\H{{\mathcal{H}}}

\def\Ka{{\mathcal{K}}}
\def\L{{\mathcal{L}}}

\def\folF{{\F}}

\def\X{{\mathcal{X}}}
\def\eX{{\fr X}}

\def\GF{\G_\folF}

\def\o{\omega}

\def\CBbb{{\mathbb C}}

\def\NBbb{{\mathbb N}}

\def\RBbb{{\mathbb R}}
\def\SBbb{{\mathbb S}}

\def\ZBbb{{\mathbb Z}}

\def\del{\partial}

\def\ov#1#2{ \frac{#1}{#2}} 
\def\parder#1#2{ \ov{\del #1}{\del #2} }
\def\dd#1#2{ \ov{d #1}{d #2} }
\def\half{\ov{1}{2}}

\def\id{\operatorname{id}}
\def\Id{\operatorname{Id}}

\def\Diff{\operatorname{Diff}}

\def\g{{\mathfrak g}}

\def\Iso{\operatorname{Iso}}

\def\proj{\operatorname{pr}}

\def\proj{\operatorname{proj}}

\def\o{\circ}
\def\X{\mathfrak X}
\def\al{\alpha}

\def\ga{\gamma}

\def\ze{\zeta}

\def\io{\iota}

\def\i{^{-1}}
\def\x{\times}

\def\Fl{\operatorname{Fl}}

\def\L{{\mathcal{L}}}

\def\g{{\mathfrak g}}

\def\pr{{\operatorname{pr}}}
\def\Id{{\operatorname{Id}}}

\def\AMSonly#1{}
\begin{document}%\topmatter
\title{
Orbifold--like and proper $\g$--manifolds
}
\author{Franz W.\ Kamber\dag\ and
Peter W.\ Michor  }
\address{
Franz W. Kamber:
Department of Mathematics,
University of Illinois,
1409 West Green Street,
Urbana, IL 61801, USA
}
\email{kamber@math.uiuc.edu }
\address{
P.\ W.\ Michor: Fakult\"at f\"ur Mathematik, Universit\"at Wien,
Oskar-Morgenstern-Platz 1, A-1090 Wien, Austria; {\it and:}
Erwin Schr\"odinger Institut f\"ur Mathematik und Physik,
Boltzmanngasse 9, A-1090 Wien, Austria
}
\email{peter.michor@esi.ac.at }
\date{July 16, 2006}
\thanks{Supported by
Erwin Schr\"odinger International Institute of Mathematical Physics,
Wien, Austria.
PWM and FWK were supported by `Fonds zur F\"orderung der wissenschaftlichen 
Forschung, Projekt P~17108~MAT'.
\\
This paper remained unfinished after the dead of Franz Kamber; now somebody wants to cite and use it
}

\begin{abstract}{
In \cite{4},\cite{5}, we generalized the concept of completion of an 
infinitesimal group action $\zeta : {\fr g} \to \X (M)$ to 
an actual group action on a (non--compact) manifold $M$~, 
originally introduced by R. Palais \cite{9}, and 
showed by examples that this completion may have quite 
pathological properties (much like the leaf space of a foliation). 
In the present paper, we introduce and investigate a tamer 
class of $\g$--manifolds, called orbifold--like, 
for which the completion has an orbifold structure.
This class of $\g$--manifolds is reasonably well--behaved with 
respect to its local topological and smooth structure to allow 
for many geometric constructions to make sense. 
In particular, we investigate proper $\fr g$--actions and generalize 
many of the usual properties of proper group actions to this more 
general setting. 
}
\end{abstract}
\def\LaTeXonly{}%\endtopmatter
%\date{\today}
%%\input amspptb.sty
%\userunningheads

\maketitle

\section{Introduction}\label{nmb:1}

\subsection{Graph foliations}\label{nmb:1.1}

Given smooth connected manifolds $B~, ~F$~, a smooth $1$--form
$\a : T (B) \to \X (F)$ defines a smooth distribution $D$ on
$B \times F \osetl{\pr_1} B$~, transversal to the fibers $F$ by lifting
tangent vectors on the base space $B$ by the formula
\begin{equation*}
D_{(x, y)} = \set{(v, \a (v) (y))}{v \in T_x (B)}~.
\end{equation*}
The distribution $D$ is integrable, that is it defines a foliation
$\folF_\a$ on $B \times F$~, if and only if $\a$ satisfies
the Cartan--Maurer equation
\begin{equation*}
d \a + \tfrac12[\a, \a] = 0~.
\tag{1}
\end{equation*}
The immersed leaves $\L \loopra  B \times F$ of $\folF_\a$~, that is
the maximal connected integral manifolds of the integrable distribution
$T \folF_\a$~, are \'etale over $B$ under composition with the
projection $\pr_1 : B\x F \to B$~.
%If $\a$ is of maximal rank, that is the evaluations
%$\ev_{x, y} (\a) : T_x (B) \to T_y (F)$ are surjective for $(x, y) \in B \x F$~,
%then the cross--sections $s : B \to B \x F~, ~s (x) = (x, y_0)$ are
%transversal to $\folF_\a$ and $\folF_\a = s^* \folF_\a$
%defines a foliation of codimension $q = \dim F$ on $B$ such that
%$T \folF_\a = \ker (\a)$~.
%At the level of microbundles, every foliation is of this type %~[Haefliger]~.

%The normal bundle $Q_{\wti{\folF}}$ of $\wti{\folF}$ is isomorphic to
%the tangent bundle $T (\pr_1)$ along the fibers of
%$F \osetl{i_x} B \x F \osetl{\pr_1} B$~:
%\begin{equation*}
%%Q_{\folF} \cong T (\pr_1) = \pr_2^* ~T (F) =
%B \x T (F) \lra B \x F~.
%\tag{2}
%\end{equation*}

If $\a : T (B) \to \X (F)$ takes values in a subalgebra
of {\it complete} vector fields, then the foliation $\folF_\a$ can be
integrated to a {\it generalized flat bundle}, that is there exists an
equivalence of bundles
\begin{equation*}
B \x F \cong \wti{B} \x_{\Ga, h_\a} F \osetl{\pi} B~, ~
\Ga = \pi_1 (B)~, ~h_\a : \Ga \to \Diff (F)
\tag{3}
\end{equation*}
such that the leaves of the foliation $\folF_\a$ correspond to the
immersions $\wti{B} \to \wti{B} \x_\Ga F$ of the form
${\ti x} \bra [{\ti x}, y]$ for $y\in F$ fixed and therefore are 
coverings of $B$~.
The homomorphism $h_\a : \Ga \to \Diff (F)$ determines the holonomy 
group $h_\a (\Ga) \subset \Diff (F)$ of the flat bundle,
respectively the foliation $\folF_\a$, and the leaves are of the form
$\wti{B} / h_\al^{-1}(\Ga_x)~,$ where $\Ga_x \subset h_\a (\Ga)$ is the
isotropy group at $x \in F$ under the action of $h_\al (\Ga)$~.
The leaf space $B\x F / \folF_\a$ of $\folF_\a$ is then given as the
orbit space $h_\al (\Ga) \bsl F~.$ It is well--known that this space
has the Hausdorff separation property, if the holonomy group $h_\al (\Ga)$
acts properly discontinuously (as a discrete group) on the fiber $F~.$
%In fact, given convergent sequences $\{ x_n \} \to x$ and
%$\{ y_n = \ga_n (x_n) \} \to y~,$ there must be a convergent, hence
%constant subsequence of $\{ \ga_n \}~,$ given by $\ga~.$ This implies
%$\{ y_n = \ga (x_n) \} \to y = \ga (x)$ and the points $x, y$ are in the
%same orbit of $h_\al (\Ga)~.$ Thus the two orbits $\bar{x} = h_\al (\Ga) (x)$
%and $\bar{y} = h_\al (\Ga) (y)$ cannot be Hausdorff separated. 
The isotropy groups $\Ga_x \subset h_\al (\Ga)$ are then finite
and the leaf space is locally Euclidean near $\bar{x}$ if
$\Ga_x = \mathbf{1}$ and has a local orbifold structure near $\bar{x}$
otherwise.

%If $\a$ is of maximal rank, the cross--sections $s$ 
%correspond to $\Ga$--equivariant submersions 
%${\ti s} : \wti{B} \to F$ onto $F$~, that is
%${\ti s} (\ti x ~\ga) = 
%h_\a (\ga)^{-1} ~{\ti s} (\ti x )~, ~\ga \in \Ga$~.
%If in addition $F$ is simply connected, the fibers of 
%$\ti s$ are connected and $\ti s$ is a {\it developing map} 
%for the foliation $\folF_\a$ on $B$~.

The main examples for the complete case occur when the fiber $F$ is
either compact, or $F = G$ is a Lie group and $\a : T (B) \to {\fr g}$
takes values in the Lie algebra $\fr g$ of left invariant vector fields
on $G$~. In the latter case the holonomy is given by a homomorphism
$h_\a : \Ga \to G$~.

\subsection{$\fr g$--manifolds}\label{nmb:1.2}

%\thetag{\thetag{\ref{nmb:1.1}.1}}
In \cite{5}, we considered a different situation, related to the integration of
infinitesimal group actions, that is the base space $B = G$ is a Lie group and
the fiber $F = M$ is a (non-compact) manifold. An infinitesimal action of $G$
on $M$ is given by a Lie homomorphism $\zeta : {\fr g} \to \X (M)$~, and
this determines a $1$--form $\a_\zeta : T (G) \to \X (M)$ satisfying the
Cartan--Maurer equation \thetag{\ref{nmb:1.1}.1} simply by composition
of $\zeta$ with the Cartan--Maurer form of $G$~.
As the examples in \cite{5} show, the completion ${_{\wti G}}M$
with respect to the simply connected Lie group $\wti G$ associated
to $\g$ may behave quite badly, e.g. ${_{\wti G}}M$ need not even
have the $T_1$--separation property.

In the present paper, we investigate a class of $\g$--manifolds
$(M, \ze)~,$ called {\it discrete} or better {\it orbifold--like},
for which the completion ${_{\wti G}}M$ has an orbifold structure.
The $\wti{G}$--completion for this class of $\g$--manifolds is reasonably 
well--behaved with respect to its local topological and smooth structure 
to allow for many geometric constructions to make sense, yet is far more 
general than the completions considered in \cite{9}.
Thus in the terminology of \cite{9}, the $\g$--manifolds considered
here need neither be univalent nor uniform.

%In section \ref{nmb:2}, we 

\section{The completion of a $\fr g$--manifold}\label{nmb:2}

\subsection{The graph of the pseudogroup}\label{nmb:2.1}

Let $M$ be a $\g$-manifold, effective and connected, so that the
action $\ze=\ze^M:\g\to \X(M)$ is injective.
Recall from \cite{1},~2.3 that the pseudogroup $\Ga(\g)$ consists of
all diffeomorphisms of the form
\begin{equation*}
\Fl^{\ze_{X_n}}_{t_n}\o\dots
     \o\Fl^{\ze_{X_2}}_{t_2}\o\Fl^{\ze_{X_1}}_{t_1}|U
\end{equation*}
where $X_i\in\g$, $t_i\in\mathbb R$, and $U\subset M$ are such that
$\Fl^{\ze_{X_1}}_{t_1}$ is defined on $U$, $\Fl^{\ze_{X_2}}_{t_2}$ is
defined on $\Fl^{\ze_{X_1}}_{t_1}(U)$, and so on.

Now we choose a connected Lie group $G$ with Lie algebra $\g$,
and we consider the integrable distribution of constant rank
($=\dim (\g)$) on $G\x M$ which is given by
\begin{equation*}
\{(L_X(g),\ze^M_X(x)):(g,x)\in G\x M, X\in\g\}\subset TG\x TM,
\tag{1}\end{equation*}
where $L_X$ is the left invariant vector field on $G$ generated by
$X\in\g$. This gives rise to a foliation $\mathcal{F}_\ze$ on $G\x M$,
which we call the {\it graph foliation} of the $\g$-manifold $M$.
Note that the graph foliation is transversal to the fibers $\{g\}\x M$ of
$\pr_1:G\x M\to G$.

Consider the following diagram, where $L_\ze (g,x) = L(g,x)$ is the leaf 
through $(g,x)$ in $G\x M$, $\mathcal{O}_\g(x)$ is the $\g$-orbit through $x$ 
in $M$, and where $W_x\subset G$ is the image of the leaf $L(e,x)$ in 
$G$. 
\begin{equation*}
\xymatrix{
 & L(e,x)  \ar[dd] \ar[dr] \ar@{->>}[rr]^{\pr_2} &  & {\mathcal{O}}_\g (x) 
\ar@{^{(}->}[d] \\
 & & G\x M \ar[d]^{\pr_1} \ar[r]^{\pr_2} & M \\
[0,1] \ar@{-->}[uur]^{\tilde c} \ar[r]^{c}  &  W_x 
\ar@{^{(}->}[r]^{\text{open}} & G & }
\tag{2}
\end{equation*}
Moreover we consider a piecewise smooth curve $c:[0,1]\to W_x$ with 
$c(0)=e$ and we assume that it is liftable to a smooth curve 
$\tilde c:[0,1]\to L(e,x)$ with $\tilde c(0)=(e,x)$. Its endpoint 
$\tilde c(1)\in L(e,x)$ does not depend on small (i.e.\ liftable to 
$L(e,x)$) homotopies of $c$ which respect the ends. This lifting 
depends smoothly on the choice of the initial point $x$ and gives 
rise to a local diffeomorphism 
$\ga_x(c): U\to \{e\}\x U\to \{c(1)\}\x U'\to U'$, 
a typical element of the pseudogroup $\Ga(\g)$ which is defined 
near $x$. 
See \cite{1},~2.3 for more information.
Note, that the leaf $L(g,x)$ through $(g,x)$ is given by 
\begin{equation*}
L(g,x)=\{(gh,y):(h,y)\in L(e,x)\}=(\mu_g\x\operatorname{Id})(L(e,x))
\tag{3}
\end{equation*}
where $\mu:G\x G\to G$ is the multiplication and 
$\mu_g(h)=gh=\mu^h(g)$. 

\subsection{The holonomy and fundamental groupoid}\label{nmb:2.2}

We briefly recall the construction of the holonomy groupoid of a foliation
$(M, \folF)$~. For $x, y \in \L$~, consider (piecewise) smooth curves
$c : x \to y$ in the leaf $\L$~, that is $\dot{c} (t) \in T(\folF)_{c (t)}$~.
By varying the initial data $(x~, \dot{c} (0))$ and considering small
variations  of $c$~, always staying in the leaves, respectively flow charts,
one constructs a local diffeomorphism $\ga : T_x \to T_y~, ~x~, ~y \in \L$, 
from  a small transversal disc $T_x \subset M$ to a small transversal disc
$T_y \subset M$~.  The germ $\ti \ga$ of $\ga$ at $x$ depends only
on the homotopy class $\{ c \}$ of the path $c$ inside the leaf,
so that we may write ${\ti \ga}_{(x, \{ c \}, y)} = [\ga]$~.
The `smooth' holonomy groupoid $\GF \to M \x M$ consists of all triples
$[x, c, y]$~, where
$c : x \to y$ is a path inside a leaf and $[x, c, y] = [x', c', y']$ iff
$x = x'~, ~y = y'$ and $c~, ~c'$ define the same holonomy transformation at
the germ level. 
%Since the leafwise homotopy relation implies
%$[x, c, y] = [x, c', y]$~, there is a canonical surjection
%$\Ga_{\folF} \to \GF$ of the homotopy or fundamental groupoid
%$\Ga_{\folF}$ to the holonomy groupoid $\GF$ of $\folF$~.
The holonomy groupoid $\GF$ is a manifold, which in general is not Hausdorff,
and the source and range maps $s, r ~:~ \GF \ra M$ determine submersions
onto the manifold $M$~, whose fibers are the holonomy coverings of the
leaves of $\folF$~.
%Intuitively, the leaf space is represented by the disjoint union of all those
%transversal discs, with the (local) holonomy transformations between them keeping
%track of the identifications which take place when leaves pass through those discs.
%A discrete version of collecting these data is called a
%Haefliger cocycle for $\folF$~.

The fundamental groupoid $\Ga_{\folF_\zeta}$ of the graph foliation
$\folF_{\zeta}$ on $G \x M$~:~
In the case of the graph foliation $\folF_{\zeta}$ on $G \x M$~,
the fundamental groupoid takes a special form, which reflects the
`incompleteness' of the infinitesimal group action.
The composition $L_{\zeta} (g, x) \looparrowright G \x M \osetl{p_1} G$
is an \'etale map over $G$~, whose range is a connected open subset
$W_{g, x} \subseteq G~, ~g \in W_{g, x}$~. By equivariance, we have
\begin{equation*}
W_{g, x} = p_1 (L_{\zeta} (g, x)) = p_1 (g L_{\zeta} (e, x)) =
g ~p_1 (L_{\zeta} (e, x)) = g W_x \subseteq G~, ~W_x = W_{e, x}~.
\end{equation*}
The projection
${\mathcal{O}}_{\fr g} (x) = {\zeta (\fr g}) (x) = p_2 ~L_{\zeta} (e, x)
\subseteq M$  is the $\fr g$--orbit through $x$~.
Since $L_{\zeta} (e, x) \osetl{p_1} W_x$ is \'etale, a path (piecewise
smooth curve) $\ti{c} : (e, x) \to (g, y)$ in $L_{\zeta} (e, x)$~,
i. e. $\dot{\ti c} (t) \in T(\folF)_{{\ti c} (t)}$~,
is uniquely determined by its projection $c = p_1 \circ \ti{c}$ in $W_x$
and the initial value $\ti{c} (0) = (e, x) \in L_{\zeta} (e, x)$~.
Such paths $c$ in $W_x$ are called {\it liftable} at $(e, x)$~.
We denote by $P_{\folF_{\zeta}, e} \subset P(G, e) \x M$ the subset of
liftable paths. This subset is open in $P(G, e) \x M$~, since
a path $c$ lifts to a unique path $\ti c$ in $L_{\zeta} (e, x)$ starting at
$(e, x)$~, provided that $c$ lies in a small convex neighborhood $U_e$
of the point path $e \in W_x$~, that is $c$ stays in a small convex
neighborhood $U_e$ of $e \in W_x \subset G$~.
Two paths $c~, ~c' \in P_{\folF_{\zeta}, e}$ with the same endpoints
are said to be related by a liftable homotopy, if their lifts
$\ti{c}~, ~\ti{c}'$ at $(e, x)$ are related by a leafwise homotopy
(covered by a finite number of flow charts) in $G \x M$~.
This implies that the lifted paths $\ti{c}~, ~\ti{c}'$ have the same
endpoint $(g, y) \in G \x M$~.
Using the $G$--invariance of the graph foliation $\folF_{\zeta}$ and its
leaves, we may translate the previous construction to any initial point
$(g, x)$ and so obtain the open subspace
$P_{\folF_{\zeta}} \subset P(G) \x M$ of liftable paths~:
\begin{equation*}
P_{\folF_{\zeta}} = \bigcup_{g \in G} ~P_{\folF_{\zeta}, g} = 
\bigcup_{g \in G} ~g ~P_{\folF_{\zeta}, e}
\end{equation*}
Identifying paths in $P_{\folF_{\zeta}}$ which are related by
liftable homotopies with fixed endpoints, we obtain exactly the
fundamental groupoid $\Ga_{\folF_{\zeta}} \ra \G_{\folF_{\zeta}}$ of
$\folF_{\zeta}$~.
By construction, the subset
$\Ga_{\folF_{\zeta}} ~(g, x) \subset \Ga_{\folF_{\zeta}}$
of classes of liftable paths starting at $(g, x)$
parametrizes the universal covering of the leaf $L_{\zeta} (g, x)$
and therefore the leaf $L_{\zeta} (g, x)$ by a covering map
\begin{gather*}
\Ga_{\folF_{\zeta}} ~(g, x) \cong \wti{L_{\zeta} (g, x)}
\osetl{\pi_{(g, x)}} L_{\zeta} (g, x)
~\qquad~,~\qquad~
\pi_{(g, x)} (\{ c \}, x) = \{ {\ti c} \} (1)~;   \\
L_{\zeta} (g, x) = \{(g', y) : (\{ c \}, x) \in \Ga_{\folF_{\zeta}} (g, x)~ ,~
\pi_{(g, x)} (\{ c \}, x) = (g', y)\}~.
\tag{1}
\end{gather*}
Moreover, there is a commutative diagram
\begin{equation*}
\xymatrix{\Ga_{\folF_{\zeta}} ~(g, x)   \ar[rr]^{{\pi_{(g, x)}}}_{}  \ar[dd]^{{{\ti p}_1}}_{}  &  & L_{\zeta} (g, x)    \ar[dd]^{{p_1}}_{} \\
 &  &   &  &  \\
\wti{W}_{g, x}  \ar[rr]^{{\ell_g \circ \pi_x \circ \ell_{g^{-1}}}}_{}  &  & W_{g, x}}
\tag{2}
\end{equation*}
where ${\ti p}_1 ((\{ c \}, x)) = [c]$ and the horizontal maps are
universal covering maps.

Since the leafwise homotopy relation implies
$[x, c, y] = [x, c', y]$~, there is a canonical surjection
$\Ga_{\folF} \to \GF$ of the 
%homotopy or 
fundamental groupoid $\Ga_{\folF}$ to the holonomy groupoid $\GF$ of $\folF$~

The pseudogroup of local diffeomorphisms on $M$ and the
corresponding smooth groupoid associated to the $\fr g$--manifold $(M, \zeta)$
are exactly the holonomy pseudogroup of local diffeomorphisms on the transversal
$M$ and the holonomy groupoid $\G_{\folF_\zeta}$~, respectively the fundamental
groupoid $\Ga_{\folF_\zeta}$ of the graph foliation $\folF_{\zeta}$
on $G \x M$~.

As is true for any foliation, the leaves $L_{\ze} (g, x)$ of
$\mathcal{F}_{\zeta}$~, viewed as immersed submanifolds of $G \times M$,
are initial submanifolds. Since leaves are by definition maximal
connected integral manifolds of the integrable distribution
$T \mathcal{F}_{\zeta}$, two leaves which intersect must be identical, that is
$L_{\zeta} (g_1, x_1) \cap L_{\zeta} (g_2, x_2) \neq \emptyset$ implies
$L_{\zeta} (g_1, x_1) = L_{\zeta} (g_2, x_2)$.

\subsection{Enlarging to group actions}\label{nmb:2.3}

In the situation of \ref{nmb:2.1} let us denote by
${_G}M = G \x_\g M = G \times M / \mathcal{F}_{\zeta}$ the
space of leaves of the foliation $\mathcal{F}_{\zeta}$ on $G\x M$,
with the quotient topology. For each $g\in G$ we consider the mapping
\begin{equation*}
j_g:M \East{{i_g}}{}\{g\}\x M \subset G\x M \East{{\pi}}{}G\x_\g M.
\tag{1}
\end{equation*}
Note that the submanifolds $\{g\}\x M \subset G\x M$ are transversal
to the graph foliation ${\mathcal{F}}_\ze$~.  
The leaf space ${_G}M$ of ${\mathcal{F}}_\ze$ admits a unique smooth structure, possibly
singular and non--Hausdorff, such that a mapping $f : {_G}M \to N$
into a smooth manifold $N$ is smooth if and only if the compositions
$f \circ j_g ~: ~M \to N$ are smooth. For example we may use the the
structure of a {\it Fr\"olicher space} or {\it smooth space} in the
sense of \cite{8},~section~23 on ${_G}M=G\x_\g M$.
%The canonical open maps $j_g ~: ~M \to {_G}M~, ~g \in G$
%are called the charts of ${_G}M$~. 
We write $j_g (x) = [g, x] \in {_G}M$~, where $[g, x]$
represents the leaf $L_{\zeta} (g, x)$
of ${\mathcal{F}}_{\ze}$ through $(g, x)$~. 
We shall also use the notation
$V_g=j_g(M)$~, an open subset of ${_G}M$.
The leaf space ${_G}M = G \times_{\g} M$ is a smooth $G$--space by
$\ell_g [g', x] = g [g', x] = [gg', x]$~,
the action being induced by the left action of $G$ on $G \times M$~.
By construction, we have
\begin{equation*}
j_{g'} = \ell_{g'g^{-1}} \circ j_{g}~,
\tag{2}
\end{equation*}
and therefore the relationships between the maps 
$j_g~, ~g \in G$ are given by left translation in ${_G}M$~.

By construction, for each $x\in M$ and for $g'g\i$ near enough to $e$
in $G$ there exists a curve $c:[0,1]\to W_x$ with $c(0)=e$ and
$c(1)=g'g\i$ and an open neighborhood $U$ of $x$ in $M$ such that
for the smooth transformation $\ga_x(c)$ in the pseudogroup
$\Ga(\g)$ we have
\begin{equation*}
j_{g'}|U = j_{g}\o \ga_x(c).
\tag{3}
\end{equation*}
Thus the mappings $j_g$ may serve as a replacement for charts in the
description of the smooth structure on ${_G}M$ and we therefore call the maps 
$j_g : M \to V_g \subset {_G}M, g \in G$ the charts of ${_G}M$~. 
Note that the mappings $j_g$ are not injective in general. Even if $g=g'$
there might be liftable smooth loops $c$ in $W_x$ such that
\thetag{3} holds. Note also some similarity of the system of
`charts' $j_g$ with the notion of an orbifold where one uses
finite groups instead of pseudogroup transformations.

\begin{thm*} {\rm See \cite{5},~theorems~5, 7, and 9.}
The $G$--completion ${_G}M$ has the following properties:

\therosteritem{4}
  Given any Hausdorff $G$-manifold $N$ and $\g$-equivariant mapping 
  $f:M\to N$ there exists a unique $G$-equivariant continuous mapping 
  $\tilde f:{_G}M\to N$ with $\tilde f\o j_e = f$. Namely, the mapping 
  $\bar f:G\x M\to N$ given by $\bar f(g,x)=g.f(x)$ is smooth and 
  factors to $\tilde f:{}_GM\to N$.

\therosteritem{5}
  If $M$ carries a symplectic or Poisson structure or a Riemannian metric 
  such that the
  $\g$--action preserves this structure or is even a Hamiltonian action
  then the structure `can be extended to ${}_GM$
  such that the enlarged $G$-action preserves these structures or is even
  Hamiltonian'.

\therosteritem{6}
  Suppose that the $\g$-action on $M$ is transitive and $M$ is connected. 
  Then there exists a (possibly non-closed) subgroup
  $H\subset G$ such that the $G$-completion ${_G}M$ is diffeomorphic to
  $G/H$. In fact, 
  $H=\{g\in G: (g,x_0)\in L(e,x_0)\} =\{g\in G: L(g,x_0)= L(e,x_0)\}$, and
  the Lie algebra of the path component of $H$ equals $\g_{x_0}=\{X\in
  \g:\ze_X(x_0)=0\}$.

\therosteritem{7}
  In general, the $G$-completion ${_G}M$ is given as follows:
  Form the leaf space $M/\g$, a quotient of $M$ which may be
  non-Hausdorff and not $T_1$ etc.
  For each point $z\in M/\g$, replace the orbit $\pi\i(z)\subset M$ by
  the homogeneous space $G/H_x$ described in 
  %loc. cit. 
  \cite{5}, theorem 7, 
  where $x$ is some point in the orbit $\pi\i(z)\subset M$, chosen 
  lo\-cal\-ly in transversals to the $\g$-orbits in $M$. 
\end{thm*}

\subsection{Intersection sets}\label{nmb:2.4}

The transversal $\{g'\}\x (M) \subset G \times M$ intersects the leaf
$\io : L_{\zeta} (g, x) \hra G \x M$ in the set
$\{g'\}\x I (g'; g, x)$~, that is
\begin{equation*}
\{g'\}\x I (g'; g, x) := (\{g'\}\x M) \cap \io L_{\zeta} (g, x)
= \io J (g'; g, x),
\end{equation*}
where $J (g'; g, x) \subset L_{\zeta} (g, x)$.
Since $L_{\zeta} (g, x) = (\mu_g\x \Id) ~L_{\zeta} (e, x)$ we have
\begin{align*}
\{g'\}\x I (g'; g, x) &= (\{g'\}\x M) \cap (\mu_g\x \Id) L_{\zeta}(e, x) \\
&= (\mu_g\x \Id) (\{g^{-1} g'\}\x M \cap L_{\zeta} (e, x) ) \\
&= (\mu_g\x \Id) (\{g^{-1} g'\}\x I (g^{-1} g'; e, x)), \\
\{g'\}\x I (g'; g, x)
     &= (\mu_{g'}\x \Id)~ (\{e\}\x M \cap L_{\zeta} (g, x)  \\
&= (\mu_{g'}\x \Id) ~(\{e\}\x M \cap L_{\zeta} ({g'}^{-1} g, x) ) \\
&= (\mu_{g'}\x \Id) ~(\{e\}\x I (e; {g'}^{-1} g, x), \\
\end{align*}
and therefore
\begin{equation*}
I (g'; g, x) = I (g^{-1} g'; e, x) = I (e; {g'}^{-1} g, x) ~\qquad~,
~\qquad~ I (e; g, x) = I (g^{-1}; e, x)~.
\tag{1}
\end{equation*}
In terms of intersection sets, the construction in \ref{nmb:2.2} may now be
written as
\begin{equation*}
y \in I (e; g, x) = I (g^{-1}; e, x) \lBia
(e, y) = \{ \ti c \} (1)~, ~(\{ c \}, x) \in \H_{\folF_{\zeta}} ~(e; g, x)~.
\tag{2}
\end{equation*}
%%Thus the intersection sets $I (e; g, x) = I (g^{-1}; e, x)$
%record exactly the possible `values' for $y = g x$ of the
%integrated `action' of $G$ on $M$~.

\subsection{Recurrence sets}\label{nmb:2.5}

By definition, the intersection sets $I (g'; g, x) \subset M$ depend
only on the leaf $L_{\zeta} (g, x)$ and the transversal $i_{g'} (M)$~,
that is we have
\begin{equation*}
I (g; g_1, x_1) = I (g; g_2, x_2)~, ~\text{for} ~[g_1, x_1] = [g_2, x_2]~.
\tag{1}
\end{equation*}
If $y \in I (g; g_1, x_1) \cap I (g; g_2, x_2)$~, then
$(g, y) \in L_{\zeta} (g_1, x_1) \cap L_{\zeta} (g_2, x_2)$ and
$[g, y] = [g_1, x_1] = [g_2, x_2]$~. Therefore distinct leaves have
disjoint
intersection sets, that is
\begin{equation*}
I (g; g_1, x_1) \cap I (g; g_2, x_2) = \emptyset~, ~
\text{for} ~[g_1, x_1] \neq [g_2, x_2]~.
\tag{2}
\end{equation*}
The following properties are equivalent~: ~
\begin{enumerate}
\item[(3).]
$[g_1, x_1] = [g_2, x_2]$~;
\item[(4).]
$L_{\zeta} (g_1, x_1) = L_{\zeta} (g_2, x_2)$~;
\item[(5).]
$g_2^{-1} g_1 ~L_{\zeta} (e, x_1) = L_{\zeta} (e, x_2)$~;
\item[(6).]
$x_1 \in I (g_1; g_2, x_2) = I (e; g_1^{-1} g_2, x_2)
~\qquad~,~\qquad~
x_2 \in I (g_2; g_1, x_1) = I (e; g_2^{-1} g_1, x_1)~.$
\end{enumerate}
Further the construction \ref{nmb:2.1} implies the transitivity relation
\begin{equation*}
\left.
\aligned &x_1 \in I (g_1; g_2, x_2) = I (e; g_1^{-1} g_2, x_2)\\
         &x_2 \in I (g_2; g_3, x_3) = I (e; g_2^{-1} g_3, x_3)
         \endaligned\right\}
~\Longrightarrow~ x_1 \in I (g_1; g_3, x_3) = I (e; g_1^{-1} g_3, x_3)~.
\tag{7}
\end{equation*}
From \thetag{\ref{nmb:2.4}.1}, we have
$x \in I (g; g, x) = I (e; e, x)~, ~\forall x \in M$. 
We call the sets $I (g; g, x)$ {\it recurrence sets}, because they
count the the points at which the leaf $L_{\zeta} (g, x)$ through
$(g,x)$ returns to the transversal $\{g\}\x M$~.
 From (7) and the previous formula, it follows that the recurrence
sets are in addition symmetric and transitive~: ~
\begin{gather*}
y \in I (g; g, x) ~\lBia~ x \in I (g; g, y) ~, \\
y \in I (g; g, x)~, ~z \in I (g; g, y) ~\lRa~ z \in I (g; g, x)~. 
\end{gather*}
Therefore the relation
\begin{equation*}
x \sim y ~\lBia~ y \in I (g; g, x) = I (e; e, x)
\tag{8}
\end{equation*}
is an equivalence relation on $M$~, independent of $g \in G$~. 

The fiber of $j_g : x \mapsto (g, x) \mapsto [g, x]$ is
exactly the recurrence set $I (g; g, x) = I (e; e, x)$~. 
%and therefore independent of $g \in G$~. 
Since the projection $G \x M \to {_G}M$
to the leaf space is an (open) identification map by definition,
we see that $j_g : M \to V_g$ is equivalent to the identification map
given by the equivalence relation \thetag{\ref{nmb:2.5}.8}.
Further it follows from \thetag{\ref{nmb:2.4}.2} that the units
$\Ga_{\folF_{\zeta}} ~(e; e, x_0)$ in the fundamental groupoid
$\Ga_{\folF}~,$ respectively the units
$\H_{\folF_{\zeta}} ~(e; e, x_0)$ in the holonomy groupoid
$\H_{\folF}~,$ act locally near $x_0 \in M~,$ with the
intersection sets $I (e; e, x)$ as orbits at $x \in M$ near $x$~.

\subsection{Uniformity of the recurrence sets}\label{nmb:2.6}

It will be desirable to have some control on 
how the recurrence sets $I (g; g, x) = I (e; e, x)$ change 
under variations of $x\in M$ in a small
transversal disc $\{g\}\x U_x~, ~U_x \subset M$~.
We set $g = e$ and represent the holonomy element
${\ti \ga}_{(x_0, \{ c \}, y_0)} = [\ga]$ by a suitable local
diffeomorphism $\ga : U_{x_0} \to U_{y_0}$~.
Then $y_0 = \ga (x_0) \in I (e; e, x_0)$ and it follows from
the holonomy construction that
$y = \ga (x) \in I (e; e, x)~, ~\forall x \in U_{x_0}$~. 
%In fact, given $y_0 \in I (e; e, x_0)$~, we have constructed a 
%continuous slice around $y_0$
%\begin{equation*}
%%\{ y = \ga (x) \mid x \in U_{x_0} \} \subset
%\bigcup_{x \in U_{x_0}} ~I (e; e, x)~,
%~y_0 = \ga (x_0)~, ~y_0 = \ga (x_0)~.
%\end{equation*}
The open neighborhood $U_{x_0} \subset M$ is not necessarily 
uniform in $y_0 = \ga (x_0) \in I (e; e, x_0)$ and other elements 
may enter or leave the recurrence set $I (e; e, x)$ as we vary $x$~.

We say that the set 
%of units 
$\H_{\folF_{\ze}} (e; e, x_0)$ in the
holonomy groupoid, respectively the
recurrence sets $I (e; e, x)$ are {\it uniform} at $x_0\in M~,$
if the elements ${\ti \ga}\in\H_{\folF_{\zeta}} ~(e; e, x_0)$
can be represented by a pseudogroup element
$\ga : U_{x_0} \oset{\cong} U_{\ga (x_0)}=\ga U_{x_0}~,$
with $U_{x_0}$ being independent of ${\ti \ga}~.$ 
If in addition $U_{x_0}$ can be chosen such that 
$\{ x_n \} \to x_0$ in $U_{x_0}$ and $\{ \ga_n (x_n) \} \to y_0$ 
implies $\ga_n = \ga$ for a (cofinal) subsequence of $\{ \ga_n \}$~, 
we say that $\H_{\folF_{\ze}} (e; e, x_0)$ acts 
{\it properly discontinuously} near $x_0\in M~.$  
It follows in particular that the recurrence sets 
$I (e; e, x) \subset M$ are discrete and closed and that the 
intersections $U_{x_0} \cap \ga (U_{x_0}) = \emptyset~,$ 
except for finitely many elements $\ga~.$ 
Therefore the isotropy group $\Ga_{x_0}$ at $x_0\in M$~, given by the units in  
$\H_{\folF_{\ze}} (e; e, x_0)$ 
is finite and we may choose $U_{x_0}$ invariant under 
$\Ga_{x_0}~.$ 
Given any $y_0 \in I (e; e, x_0)$~, we have now constructed a 
uniform continuous slice (or section) around $y_0$~, that is
\begin{equation*}
\{ y = \ga (x) \mid x \in U_{x_0}, \ti\ga \in \H_{\folF_{\ze}}~(e; e, x_0 \} 
= \bigcup_{x \in U_{x_0}} ~I (e; e, x)~.
\tag{1}
\end{equation*}
Equivalently, there is a saturated neighborhood $U$ of
$L (e, x_0) \subset G\x M$ of the form
\begin{equation*}
U \cong \what{L} (e, x_0) \x_{\Ga_{x_0}} U_{x_0} \to 
\Ga_{x_0} \bsl U_{x_0}~,
\tag{2}
\end{equation*}
where $\what{L} (e, x_0) \to L (e, x_0)$ is a finite 
covering space of $L (e, x_0)$ with group $\Ga_{x_0}~,$ 
that is the saturated neighborhood $U$ of the leaf $L (e, x_0)$ 
admits a Seifert fibration over $\Ga_{x_0} \bsl U_{x_0}$~. 
If the above properties are valid for any $x_0 \in M$~, we 
use the terms without reference to the basepoint $x_0$~. 

As $[e, x_0]\in V_e$ approaches the boundary
$\del V_e = \wbar{V}_e\setminus V_e$ of $V_e \subset {_G}M~,$
the intersection sets may change drastically. 
This leads to the notion of limit elements explained in
\ref{nmb:2.8} below.

\subsection{Hausdorff separation property for ${_G}M$}\label{nmb:2.7}

Recall that $V_g=j_g(M)\subseteq {_G}M$.
Suppose that the elements $z = [g, x] \in V_g~, ~z' = [g', y] \in V_{g'}~, ~z \neq z'$
are not Hausdorff separated. Then both $z$ and $z'$ would have to be arbitrarily
close to points $z_n \in V_g \cap V_{g'}$~, that is leaves which intersect both
transversals $\{g\}\x M$ and $\{g'\}\x M$ in $G\x M$, and we would have
$z_n = [g, x_n] = [g', y_n] \in V_{g} \cap V_{g'}$~, for convergent
sequences $\{ x_n \} \to x~, ~\{ y_n \} \to y$~.
Therefore the Hausdorff property for the leaf space ${_G}M$ is equivalent
to the following statement~:~

Suppose that $\{ x_n \} \to x$ and $\{ y_n \} \to y$ are convergent
sequences in $M$, such that $(g', y_n) \in L_{\zeta} (g, x_n)$. Then we have
$(g', y) \in L_{\zeta} (g, x)$. Equivalently, $y_n \in I (g'; g, x_n)$
implies $y \in I (g'; g, x)$. 

The Hausdorff separation property for $V_e \subseteq {_G}M,$
and hence for any $V_g \subseteq {_G}M,$
is equivalent to the following statement~:

Suppose that $\{ x_n \} \to x$ and $\{ y_n \} \to y$ are convergent
sequences in $M$, such that $(e, y_n) \in L_{\zeta} (e, x_n)$. Then we have
$(e, y) \in L_{\zeta} (e, x).$ Equivalently, $y_n \in I (e; e, x_n)$
implies $y \in I (e; e, x).$

\subsection{Limit elements}\label{nmb:2.8}

The incompleteness of the $\fr g$--action must of course show up somewhere
and this happens if we move towards the boundary
$\del V_{g} = \wbar{V}_g \bsl V_g$ of the open sets $V_g \subseteq {_G}M$~.
We say that an element $\hat z  = [{\hat g}, {\hat x}] \in {_G}M$ is a
{\it limit element} for the transversal $i_{g_0} (M)$ if
$\hat z \in \del V_{g_0}$~.
This means that the leaf $L_{\zeta} ({\hat g}, {\hat x})$~,
representing $\hat z$ does not intersect the transversal $i_{g_0} (M)$~,
that is $I (g_0; {\hat g}, {\hat x}) = \emptyset$~, but is arbitrarily close
to a leaf $L_{\zeta} (g_n, x_n)$ such that $I (g_0; g_n, x_n) \neq \emptyset$~.
%Equivalently, we may require that $L_{\zeta} ({\hat g}, {\hat x})$
%intersects transversals $i_{g_n} (M)$~, which are arbitrarily close to
%$i_{g_0} (M)$~, that is $I (g_n; {\hat g}, {\hat x}) \neq \emptyset$~,
%but $I (g_0; {\hat g}, {\hat x}) = \emptyset$~.
It is obviously sufficient to elaborate this property for $g_0 = e$
and $V_e \cap V_{\hat g} \neq \emptyset$~.

The following statements are equivalent~: ~
\begin{enumerate}
\item[(1).]
     $\hat z = [{\hat g}, {\hat x}]  \in {_G}M$ is a limit element
     for the transversal $i_e (M),$ that is
     $\hat z \in \del V_e \cap V_{\hat g}~;$
\item[(2).]
     There exist sequences $\{ x_n \}~, ~\{ y_n \}$ in $M$~,
     such that $\{ x_n \}$ converges to ${\hat x}$ and
     $y_n \in I (e; {\hat g}, x_n) = I (e; e, y_n)$~, but
     $I (e; {\hat g}, {\hat x}) = \emptyset$~
\end{enumerate}
%The Hausdorff property for $V_e$ implies that the sequence
%$\{ y_n \}$ in $M$ cannot have an accumulation point for any choice of
%the $y_n \in I (e; {\hat g}, x_n)$~. Indeed, if there were a convergent sequence
%$\{ y_n \} \to y$~, such that $y_n \in I (e; {\hat g}, x_n)$~, we would have
%$[{\hat g}, {\hat x}] = [e, y] \in V_e$ and $[{\hat g}, {\hat x}]$ would not
%be a limit element for $V_e$~.

\subsection{Separation of limit elements}\label{nmb:2.9}

If $V_{\hat g} \cap V_{e} = \emptyset$~, then any points
$z_0 \in V_{e}$ and ${\hat z} \in V_{\hat g}$ are obviously
Hausdorff separated.
Suppose that $V_e \cap V_{\hat g} \neq \emptyset~, ~{\hat g} \neq e$~.
If $z_0 \notin {\wbar V}_{\hat g}$~, then $\hat z$ and $z_0$
are still Hausdorff separated by the disjoint open sets $V_{\hat g}$ and
$V_{e} \setminus (V_{e} \cap {\wbar V}_{\hat g}) \subset V_{e}$~.
The same is true if $\hat z \notin {\wbar V}_{e}$~.
For $z_0~, ~{\hat z} \notin V_e \cap V_{\hat g}$~,
the Hausdorff property can therefore fail only for
limit elements $z_0 \in \del V_{\hat g}$ of the form
$z_0  = [e, y_0]$~, that is $z_0 \in V_e \cap \del V_{\hat g}~,$
and $\hat z \in \del V_{e}$ of the form
$\hat z = [{\hat g}, {\hat x}]$~, that is
$\hat z \in \del V_{e} \cap V_{\hat g}$~.
Then ${\hat z} \neq z_0$~,
since $(V_{e} \cap \del V_{\hat g}) \cap (\del V_{e} \cap V_{\hat g}) =
(V_{e} \cap \del V_{e}) \cap (V_{\hat g} \cap \del V_{\hat g}) = \emptyset$~.
%Under the Hausdorff property, limit elements are unique in the following sense~:
%$z_0~, ~\hat z$ cannot both be the limit of a convergent sequence
%$z_n = [e, y_n] = [{\hat g}, x_n] \in V_e \cap V_{\hat g}$ of leaves
%which intersect both transversals $i_{e} (M)$ and $i_{\hat g} (M)$~,
%that is $y_n \in I (e; {\hat g}, x_n)~, ~x_n \in I ({\hat g}; e, y_n)$ and
%$I (e; {\hat g}, {\hat x}) = I ({\hat g}; e, y_0) = \emptyset$~.
The set $H \subset {_G}M \x {_G}M \bsl \Delta ~{_G}M$ of pairs
$(z_0, {\hat z})$ of distinct non--separable points is obviously
invariant under the diagonal action of $G$ and we may write
$H = \cup_{g \in G} ~g H_e$~, where pairs $(z_0, {\hat z})$ in $H_e$ 
satisfy $z_0\in V_e~.$
The preceding argument shows that the set $H_e$ is restricted by
\begin{equation*}
H_e \subseteq V_e \x V_e \cup \left\{
\bigcup_{g \in G}^{V_e \cap V_g \neq \emptyset} ~
V_e \cap \del V_g \x \del V_e \cap V_g \right\}~.
\tag{1}
\end{equation*}
We say that pairs of limit points are {\it Hausdorff separable}
if the situation described above does not occur, in which case
we have $H_e \subseteq V_e \x V_e~.$

Non--separable limit points typically occur when a complete
$\g$--action on $M$ is made incomplete on $M_1 = M \bsl A$ by
removing a closed subset $A \subset M$ from $M$, such that
some orbits of the complete $\g$--action on $M$ become
disconnected on $M_1$. 
Such an orbit gives rise to multiple orbits in the completion
${_{\wti G}}M_1$ and non--separable pairs appear whenever there are
nearby orbits in $M$ which are still connected in $M_1 = M \bsl A$
(compare the examples in \cite{4}). 
%and example \ref{nmb:3.6} below).
Then the canonical map ${_{\wti G}}M_1 \to {_{\wti G}}M \cong M$
defines $M$ as a canonical Hausdorff quotient of ${_{\wti G}}M_1~.$

\begin{prop}\label{nmb:2.10}
The Hausdorff separation property for the leaf space
${_G}M$ is characterized as follows$~:~$

\begin{enumerate}
\item[(1).]
The leaf space ${_G}M$ is Hausdorff, if and only if
the open set $V_e \subseteq {_G}M $ is Hausdorff and
pairs of limit points are Hausdorff separable.
Then also each completion of a $\g$-orbit in $M$ in the sense of
\thetag{\ref{nmb:2.3}.7} is of the form $G/H$ for a closed subgroup $H$ of $G$.

\item[(2).]
The open set $j_e(M)=V_e \subseteq {_G}M $ is Hausdorff, if 
%the recurrence sets $I (e; e, x) \subset M$ are 
%discrete and closed and 
the action of the holonomy ${\mathcal{H}}_{\folF_{\zeta}} (e;e, x)$
at $x$ is properly discontinuous for all $x \in M~.$
\end{enumerate}
\end{prop}

\begin{demo}{Proof}
%\thetag{\thetag{\ref{nmb:2.3}.2}}
From \thetag{\ref{nmb:2.3}.2}, we see that the open sets $V_g~, ~g \in G$
have the Hausdorff property if and only if $V_e$ has the Hausdorff
property. \thetag{1} follows then from the preceding 
analysis in \ref{nmb:2.9}. To prove \thetag{2}~, let $\{ x_n \} \to x$ and $\{ y_n \} \to y$ be
convergent sequences in $M$~, such that $y_n \in I (e; e, x_n)$~.
%Since $\{ x_n \} \to x$ and $y_n \in I (e; e, x_n)$~,
%we can use the slice construction in 
%\thetag{\thetag{\ref{nmb:2.6}.1}} 
%to find open neighborhoods $V_n (x) \subset M$ shrinking to $x$ and 
%holonomy elements  $\ga_n : V_n (x) \to W_{w_n}~, ~w_n \in I (e; e, x)$~, 
%such that $x_n \in V_n (x)~, ~\ga_n (x_n) = y_n~, ~\ga_n (x) = w_n$ 
%for almost all $n$~. Since the intersection sets $I (e; e, x)$ are 
%discrete and closed by assumption, the sequence $\{ w_n \}$ has no 
%accumulation point and we must have $w_n = w$ for almost all $n$~. 
Since $\{ x_n \} \to x$ and $y_n \in I (e; e, x_n)$~,
we can use the slice construction around $x$ in \thetag{\ref{nmb:2.6}.1} 
to find an open neighborhood $U_x$ and a sequence $\ga_n$ so that 
$x_n \in U (x)~, ~\ga_n (x_n) = y_n$~. 
By assumption, the action of the holonomy groupoid is
properly discontinuous near $x$, so we have $\ga_n = \ga$ 
for a subsequence of $\{ \ga_n \}$~.
But then $y_n = \ga (x_n)$ and this converges to $y = \ga (x)$~.
Thus we have $y \in I (e; e, x)$~.
\qed\end{demo}

\subsection{Orbifold structure of $j_e : M \to V_e \subset{_G}M$}\label{nmb:2.11}
For the chart $j_e : M \to V_e$ 
%\subset {_G}M$ 
to be \'etale, that is a local diffeomorphism with respect to the 
differentiable structure induced in ${_G}M$ (cf. \ref{nmb:2.4}), 
it is necessary that the recurrence sets $I (e; e, x) \subset M~, ~x \in M$ 
are discrete, closed and uniform.
If the recurrence set $I (e; e, x) \subset M$ has an accumulation point
$y_0 \in M$~, it is clear that $j_e : M \to V_e$ cannot be a local 
homeomorphism at $y_0$~. 
A sufficient condition for the \'etale structure of $j_e : M \to V_e$
is that in addition the action of the holonomy 
${\mathcal{H}}_{\folF_{\zeta}} ~(e;e, x)$ on $M$ at $x$ is properly
discontinuous, without fixed points, for all $x \in M.$
In the sequel, we will allow fixed points, which are necessarily of finite
order, that is we will allow an {\it orbifold} structure on $V_e~,$ 
so that the description in \ref{nmb:2.6} and in particular \thetag{\ref{nmb:2.6}.1}, 
\thetag{\ref{nmb:2.6}.2} apply.
By translation invariance \thetag{\ref{nmb:2.3}.2}, %\thetag{\thetag{\ref{nmb:2.3}.2}, 
the same conditions will then hold for all charts 
$j_g : V_g \subseteq {_{G}}M = G \x_{\fr g} M~, ~
g \in G$~. 

\section{Orbifold--like $\g$--manifolds}\label{nmb:3}

A priori, the local orbifold structure described in \ref{nmb:2.11} 
depends on the choice of a Lie group $G$ with Lie algebra ${\fr g}$~, 
that is on the choice of a quotient group $\wti{G} \to G$ of the simply 
connected Lie group $\wti{G}$ associated to ${\fr g}$~. 

\subsection{Definition}\label{nmb:3.1}~
The $\fr g$--action on $M$ is {\it discretely valued} with respect to $G$,
and $(M, \ze)$ a {\it orbifold--like} $\g$--manifold with respect to $G$, 
if the chart $j_e : M \to V_e \subset{_G}M$ defines an orbifold structure 
on $j_e : M \to V_e$ as described in \ref{nmb:2.11}. 
This means that the action of the holonomy ${\mathcal{H}}_{\folF_{\zeta}} ~(e;e, x)$
on $M$ at $x$ is properly discontinuous near $x$ for all $x \in M~.$

We collect some immediate properties of orbifold--like $\g$--manifolds  
in the Lemma below. 

\begin{lem}\label{nmb:3.2}~
Let $(M, \ze)$ be an orbifold--like $\g$--manifold with respect to $G$~. 
Then the following statements hold~: ~

\begin{enumerate}

\item[(1)]
The $G$--completion $_{G}M = {G} \x_{\fr g} M$  
of $(M, \ze)$ is a smooth orbifold, such that the orbifold 
structure is determined by the charts 
$j_g : M \to V_g \subseteq {G}M$~.
More precisely, if $U_x \subset M$ is a uniform neighborhood as in
\ref{nmb:2.6}, then $j_e (U_x) \cong \Ga_x \bsl U_x~,$ where $\Ga_x$ is
the finite isotropy group at $x \in M~.$

\item[(2)]
In the absence of fixed points at $x~,$ the charts 
$j_g : M \to V_g \subset {_G}M$ 
are \'etale near $x$ and therefore local diffeomorphisms.
Therefore ${_G}M$ carries a canonical structure of a smooth manifold. 

\item[(3)]
The open orbifolds $V_g \subset {_G}M$ are always Hausdorff.  
%by Proposition \ref{nmb:2.10}. 
Moreover, the Hausdorff separation property may fail only for 
pairs of limit points.

\item[(4)]
If $i_e : M \to G \x M$ is a complete transversal of $\folF$, 
then the chart $j_e : M \to {_G}M$ is a map onto
the Hausdorff orbifold ${_G}M$~. 

\item[(5)]
All intersection sets $I (g'; g, x) \subset M$ are discrete, 
closed, uniform and at most countably infinite.

\item[(6)]
The recurrence sets $I (g; g, x) \subset M$~, that is the fibers
of the charts $j_g : M \to V_g \subseteq {_{\wti G}}M$~, 
are discrete, closed, uniform and at most countably infinite.

%\item[(7)]
%The leaf space ${_G}M = {G} \times_{\fr g} M$ satisfies the
%$T_1$--separation property.

\item[(7)]
The immersed leaves
$\iota : L_{\zeta} (g, x) \loopra G \times M$
of $\folF_{\zeta}$ are closed and thus are submanifolds.

\item[(8)]
The charts $j_g : M \to {_G}M$ are morphisms of $\fr g$--orbifolds,
that is the vector fields $\zeta (X)$ on $M$ and
${\ti \zeta} (X)$ on ${_G}M~, ~X \in {\fr g}$ are $j_g$--related.

\item[(10)]
The orbit spaces $X = {\fr g} ~\bsl~ M$ and $G ~\bsl~ {_G}M$ are
homeomorphic, with the homeomorphism
\begin{equation*}
{\bar j}~: ~{\fr g} ~\bsl~ M \cong G ~\bsl~ {_G}M
\end{equation*}
induced by the chart $j_e ~: ~M \to {_G}M$~.

\end{enumerate}
\end{lem}

\begin{demo}{Proof}
The above statements follow from the constructions in 
\ref{nmb:2.6} to \ref{nmb:2.11}. 
\qed
\end{demo}

We note that there exist $\fr g$--actions which are not 
discretely valued (cf. \cite{5}, Example 8).
There exist discretely valued $\fr g$--actions,
such that $j_e : M \to {_G}M$ is not injective, with
the fibers $I (e; e, x) \subset M$ being finite ot countably infinite
(cf. \cite{5}, Examples 4).
There are examples of discretely valued $\fr g$--actions, for which
${_{\wti G}}M$ is a non--Hausdorff manifold. The Hausdorff quotient
of ${_{\wti G}}M$ may or may not be a smooth manifold (cf. \cite{4}).

\subsection{The $\wti{G}$--completion of a ${\fr g}$--manifold}\label{nmb:3.3}

We say that ${_{\wti G}}M$ is {\it the completion} of the
$\fr g$--manifold $(M, \ze)~,$ if the $\fr g$--action is
discretely valued, that is $(M, \ze)$ is orbifold--like, 
with respect to $\wti G$. 
Then ${_{\wti G}}M$ is a smooth orbifold,
possibly with a limited failure of the Hausdorff separation
property for pairs of limit points, as described in \ref{nmb:2.9}
and Proposition \ref{nmb:2.10}.1. In this case, 
we omit the reference to $\wti{G}$ in Definition \ref{nmb:3.2}. 

Recall that $G = \wti{G} / Z$~, where $Z = \pi_1 (G) \subset \wti{G}$ is a
discrete central subgroup of the simply connected Lie group $\wti{G}$
with Lie algebra $\fr g$~.
Almost from the definition of ${_G}M = G \x_{\fr g} M$~, it follows  
that we have a canonical equivalence of identification spaces
\begin{equation*}
{_G}M = G \x_{\fr g} M \cong Z ~\bsl~ {_{\wti G}}M~,
\tag{1}
\end{equation*}
that is ${_G}M$ is the orbit space of the left action of $Z$ on the
$\wti G$--space ${_{\wti G}}M$~.
As is true for coverings, the units in the fundamental  groupoid
are related by an exact sequence
\begin{equation*}
0 \ra \Ga_{\wti{\folF}_{\zeta}} ~({\ti e}; {\ti e}, x) \lra
\Ga_{\folF_{\zeta}} ~(e; e, x) \lra Z \ra 0~,
\end{equation*}
which results in a corresponding relation between the
$\Ga_{\folF_\zeta}$--orbits,
namely the intersection sets $I_{\wti G} ({\ti e}; {\ti e}, x)$ and
\begin{equation*}
I_G (e; e, x) = \bigcup_{z \in Z} ~I_{\wti G} ({\ti e}; z, x)~.
\tag{2}
\end{equation*}

\begin{prop}\label{nmb:3.4}~

\begin{enumerate}

\item[(1)]
If ${_{\wti G}}M$ is an orbifold,
then ${_G}M = Z ~\bsl~ {_{\wti G}}M$ is an orbifold,
if and only if $Z = \pi_1 (G)$ acts properly
discontinuously on ${_{\wti G}}M$~.

\item[(2)]
If ${_G}M$ is an orbifold,
then ${_{\wti G}}M$ is an orbifold. 

\item[(3)]
If ${_{\wti G}}M$ is a manifold,
then ${_G}M = Z ~\bsl~ {_{\wti G}}M$ is a
manifold,
if and only if $Z = \pi_1 (G)$ acts properly
discontinuously and without fixed points on ${_{\wti G}}M$~. 

\item[(4)]
If ${_G}M = Z ~\bsl~ {_{\wti G}}M$ is a manifold, then
${_{\wti G}}M$ is a manifold and the canonical projection
$\pi : {_{\wti G}}M \to {_G}M \cong Z ~\bsl~ {_{\wti G}}M$ 
is a covering map.

\item[(5)]
If ${_G}M = Z ~\bsl~ {_{\wti G}}M$ is Hausdorff, then ${_{\wti G}}M$
is Hausdorff. 

\end{enumerate}
\end{prop}

\begin{demo}{Proof}
\thetag{1} to \thetag{3} follow from \thetag{\ref{nmb:3.3}.1} and 
\thetag{\ref{nmb:3.3}.2}. 
\thetag{4} follows from the fact that the charts 
$j_g : M \to {_{\wti G}}M$ and $j_{Zg} : M \to {_G}M$ 
are related by $j_{Zg} = \pi \circ j_g$~. 
\thetag{5} follows from the fact that the canonical projection
$\pi$ is an identification map and hence open. 
\qed
\end{demo}

Thus the finite valuation condition, respectively the orbifold--like 
condition is the weakest condition to impose for a ${\fr g}$--manifold.   
From the preceding construction, it is easy to produce examples where 
${_{\wti G}}M$ is a manifold, but ${_G}M$ is not Hausdorff and not
locally euclidean.

\subsection{Complete $\fr g$--actions}\label{nmb:3.5}

If the $\fr g$--action $\zeta : {\fr g} \to \eX (M)$
is given by {\it complete} vector fields on $M$~, then the $\fr g$--action
integrates to a $\wti G$--action on $M$, see \cite{9}.
More precisely, we have then $W_x = G~, ~\forall x \in M~,$ the projection
$\pr_1$ in diagram \thetag{\ref{nmb:2.1}.2} is a universal covering map
and $\folF_{\zeta}$ determines a generalized flat bundle
\begin{equation*}
G \x M \ocong{\bar\psi} \wti{G} \x_Z M \osetl{\pr_1} G~, ~
Z = \pi_1 (G) \subset Z (\wti G)~,
\tag{1}
\end{equation*}
as a fiber space over $G$~. Here $\bar\psi (g, x) = [g, g x]$ is the
diffeomorphism induced by ${\ti \psi} : \wti{G} \x M \cong \wti{G} \x M~, ~
{\ti \psi} (g, x) = (g, g x)$~. Note that $\ti \psi$ transforms the right
diagonal action $(g', x) \cdot g = (g' g, g^{-1} x)$ into the right action
$(g', x) g = (g' g, x)$ and the left action $g (g', x) = (g g', x)$ into
the left diagonal action $g \cdot (g', x) = (g g', g x)$~.
Then we have from \thetag{\ref{nmb:3.3}.1}
\begin{equation*}
{_G}M = G \x_{\fr g} M \cong G \x_{\wti G} M \ocong{\psi} Z ~\bsl~ M
\tag{2}
\end{equation*}
as $G$--spaces. The recurrence sets $I (g; g, x) = I (e; e, x) = \{ Z x \}$
are exactly the $Z$--orbits of $x$ and ${_G}M$ is a manifold
if and only if $Z = \pi_1 (G)$ acts properly discontinuously and without
fixed points on $M$~.
The foliation $\folF_{\zeta}$ is induced by the submersion
$\pr_1 \circ \ti\psi : \wti{G} \x M \to \wti{G}$~,
and the $G$--action on ${_G}M$ is induced by the $\wti G$--action on $M$~,
that is by the composition
${\ti \mu} = \pr_2 \circ \ti\psi : \wti{G} \x M \to M$~.
For $G = \wti{G}$~, the flat bundle \thetag{\ref{nmb:2.1}.2} is a product
and \thetag{1}, \thetag{2} take the form
\begin{equation*}
\wti{G} \x M \ocong{\bar\psi} \wti{G} \x M \osetl{\pr_1} \wti{G}
~\qquad~,~\qquad~
{_{\wti G}}M \ocong{\psi} M~.
\tag{3}
\end{equation*}
in agreement with \thetag{\ref{nmb:3.3}.1}. Then the intersection sets $I (e; g, x) = I (g^{-1}; e, x) = \{ y \}$ are
singleton sets and the $\wti G$--action is uniquely determined by $g x = y$~,
that is $(e, y) \in L_{\zeta} (g, x)$~.
The preceding statements hold in particular if $M$ is compact
and $\eX_c (M) = \eX (M)$~.

\subsection{Example}\label{nmb:3.6}

The following example illustrates how orbifolds occur naturally
during the completion process.

Let $\g=\RBbb$ with basis $X = \dd{\ }{t}~,$ let
$M_2=\{z\in \CBbb : \abs{z} < 2 \}$ be the open disc of radius $2~,$
and let $\zeta : \g \to \X(M)$ be given by the complete vector field
\begin{equation*}
\zeta_X (z)= \ov{2\pi \iota}{n} ~z ~\parder{\ }{\theta}~, ~
z=r ~e^{2\pi\iota \theta}~, ~n \in \NBbb~, ~n>2~,
\tag{1}
\end{equation*}
with integral curves
$
z (t) = z_0 ~e^{\ov{2\pi\iota}{n} ~t}~, ~
z_0=r_0 ~e^{2\pi\iota \theta_0}~,
~z_0 \in M_2~.$
As $\ze_X$ is complete on $M_2~,$ we have ${_\RBbb}M_2 \cong M_2$
and ${_{\SBbb^1}}M_2 \cong \ZBbb \bsl M_2$ from \ref{nmb:3.5}. 
For $\SBbb^1 \cong \RBbb / \ZBbb \subset \CBbb$ the unit circle,
$\SBbb_{(n)}^1 \cong \RBbb / n \ZBbb \to \SBbb^1$ the $n$--fold
covering determined by $z \bra z^n$ and 
$\Ga = \ZBbb_n = \ZBbb / n\ZBbb$
the cyclic group of order $n>2~,$ we have moreover from
\thetag{\ref{nmb:3.5}.1}
\begin{equation*}
\SBbb^1 \x M_2 \cong \RBbb \x_{\ZBbb} M_2 \cong
\SBbb_{(n)}^1 \x_\Ga M_2 \to \SBbb^1~.
\tag{2}
\end{equation*}
It follows that the foliation $\folF_\ze$ on $\SBbb^1 \x M_2$
is determined by the flow lines in the mapping torus of
$f = e^{\ov{2\pi\iota}{n}}: M_2 \to M_2$ and that the completion
${_{\SBbb^1}}M_2 \cong \Ga \bsl M_2$ is the standard orbifold.
The leaves of the graph foliation on $\SBbb^1 \x M_2$ are closed,
parametrized by $\SBbb_{(n)}^1$ away from the origin in $M_2~,$
while the leaf at $0 \in M_2$ is parametrized by $\SBbb^1~.$
The orbit space $M_2 / \g$ is given by the half open interval
$[0, 2)~.$

To make $\ze_X$ incomplete, we remove the relatively closed annular
wedge $A = \{ (r,\theta): 1\leq r< 2, 0\leq\abs{\theta}\leq\ov{1}{n}\}$
of width $\ov{2}{n}$ around the $x$--axis from the open disc $M_2~,$
that is we set $M=M_2 \bsl A~.$
%the interval $I = \{ (x,0): 1\leq x < 2 \}$ from the open disc
%$M_2~,$ that is we set $M=M_2 \bsl I~.$
For $0 \leq r_0 < 1~,$ the integral curves (orbits)
are defined for all values of the parameter $t~,$ while
for $1 \leq r_0 < 2~,$ the integral curves are defined only for
$1<\abs{t+ n \theta_0} < n-1~, z_0 \in M~.$
Since $\SBbb^1 \cong \RBbb / \ZBbb$ via the exponential map, we have
${_{\SBbb^1}}M \cong \ZBbb \bsl {_\RBbb}M$ from \thetag{\ref{nmb:3.3}.1},
and it is sufficient to determine the $\RBbb$--completion ${_\RBbb}M$
with the induced $\RBbb$--action.
By construction of the graph foliation ${\mathcal{F}}_{\ze}$ in
\thetag{\ref{nmb:2.1}.2}, the leaves of ${\mathcal{F}}_{\ze}$ on $\RBbb \x M$
are determined explicitly as follows.
For a regular curve $\vphi = c (t)$ in $\RBbb$ starting at $\vphi_0 = c (0)~,$
we have $\dot{\vphi} (t) = \dot{c} (t) ~X~, ~\dot{c} (t)>0$ and the lifted
curve $(c (t), z (t))$ is in the leaf $L (c_0, z_0)~, ~z_0 \in M~,$
if and only if it satisfies the first order ODE
$
(z (t), \dot{z} (t)) = \dot{\vphi} (t) ~\zeta_X (z (t))
$
with initial value $z (0) = z_0 = (r_0, \theta_0) \in M$.
Substituting \thetag{1} into this equation,
we obtain the linear equation
$\dot{z} (t) =
\ov{2\pi \iota}{n} ~\dot{\vphi} (t) ~z (t)$
and therefore
\begin{equation*}
z (\vphi) = z_0 ~e^{\ov{2\pi\iota}{n} (\vphi (t) - \vphi_0)} =
r_0 ~e^{\ov{2\pi\iota}{n} (\vphi - (\vphi_0 - n \theta_0))}~.
\tag{3}
\end{equation*}
This is independent of the parametrization $\vphi = c (t)~,$
except for initial values, so we may take $\vphi = t+\vphi_0~.$
For $0\leq r_0<1~,$ the projection
\begin{equation*}
\pr_1:L (\vphi_0, z_0) \to W_{(\vphi_0,z_0)} = \RBbb~, ~0 \leq r_0 < 1
\tag{4}
\end{equation*}
is a diffeomorphism for all values of $\vphi_0~.$

The orbits of the $\g$-action are determined by the leaf structure
via $\pr_2$ in diagram \thetag{\ref{nmb:2.1}.2} and they look as follows:
The point $0 \in M \subset \CBbb$ is a fixed point. For $0<r_0<1~,$
the orbits are circles, parametrized by $\SBbb_{(n)}^1 \cong \RBbb / n \ZBbb$
as in the complete case and the orbit space of the open disc
$M_1 = \{ z\in \CBbb : \abs{z}<1 \} \subset M$ is the half open
interval $[0, 1)~.$ In other words, the foliation $\folF_\ze$ on
$\SBbb^1 \x M_1$ is again given by the mapping torus of
$f = e^{\ov{2\pi\iota}{n}}: M_1 \to M_1$ as in \thetag{2},
the $\RBbb$--completion ${_{\RBbb}}M_1 \cong M_1$ and the
$\SBbb^1$--completion ${_{\SBbb^1}}M_1 \cong \Ga \bsl M_1$
is again the standard orbifold.

Incompleteness occurs whenever the integral curve $z (\vphi)$
approaches the wedge $A \subset M_2$ in finite time, that is
$z (\vphi) \to (r_0, \pm \ov{1}{n} )~, ~1 \leq r_0 < 2~,$ or equivalently
$\vphi - (\vphi_0 - n \theta_0) \da 1~,$ respectively
$\vphi - (\vphi_0 - n \theta_0) \ua n-1~.$
It follows that the leaf $L (\vphi_0, z_0)$
is parametrized by $\vphi - (\vphi_0 - n \theta_0) \in (1, n-1)$ and that
\begin{equation*}
\pr_1:L (\vphi_0, z_0) \to W_{(\vphi_0, z_0)} = \vphi_0 + W_{z_0}~, ~
W_{z_0} = n (\ov{1}{n} - \theta_0, 1-\ov{1}{n} - \theta_0) \subset \RBbb~,
\tag{5}
\end{equation*}
given by $(\vphi, z(\vphi)) \bra \vphi~,$ is a
diffeomorphism for $1 \leq r_0 < 2~, 1 < n \theta_0 < n-1~.$
The intervals $W_{z_0}$ of length $n-2$ range from $(0, n - 2)$ to
$(2 - n, 0)~,$ with $(1-\ov{n}{2}, \ov{n}{2}-1)$ for $\theta_0 = \half~.$
This reflects on the fact that the cut--open circles are intervals and
therefore simply connected; the completion process simply adds in
copies of $M,$ rotated by the parameter $\vphi~.$
For $1\leq r_0<2~,$ we pa\-ra\-met\-rize the space of leaves
${_{\RBbb}}(M \bsl M_1)$ by eliminating the parameter
$\theta_0~.$ From \thetag{5}, we see that
$
L (\vphi_0', z_0') = L (\vphi_0, z_0),
$
if and only if $r_0' = r_0~, ~\vphi_0' - n \theta_0' = \vphi_0 - n \theta_0~.$
Setting $\theta_0' = \half~,$
%$\vphi_0 = n (\theta_0 - \ov{1}{2})$ or $\theta_0 = \ov{1}{2} + \ov{\vphi_0}{n}$
it follows from \thetag{\ref{nmb:2.1}.3} that
\begin{equation*}
L (\vphi_0, z_0) = L (\vphi_0'+\ov{n}{2}, (r_0, \half)) =
\vphi_0' + L (\ov{n}{2}, (r_0, \half))~, ~\vphi_0' = \vphi_0 - n \theta_0~.
\tag{6}
\end{equation*}
Therefore the leaves of the form
$L (\vphi_0 +\ov{n}{2}, (r_0, \half)) = \vphi_0 + L (\ov{n}{2}, (r_0, \half))$
are distinct for different values of $\vphi_0 \in\RBbb$ and
$r_0 \in [1,2)~,$ that is $\wti{G} = {\mathbb R}$ acts without
isotropy on ${_{\RBbb}}(M \bsl M_1)$.
From \thetag{6},
we see that the $\RBbb$-completion ${_{\RBbb}}(M \bsl M_1)$ has a
section over the orbit space ${_{\RBbb}}(M \bsl M_1) / \g~,$ given by
$r_0 \mapsto L (\ov{n}{2}, (r_0, \half))~.$ 
Therefore
${_{\RBbb}}(M \bsl M_1) \cong \RBbb \x ({_{\RBbb}}(M \bsl M_1) / \g)
\cong \RBbb \x [1, 2)$ and
${_{{\SBbb}^1}}(M \bsl M_1) \cong \SBbb_{(n)}^1 \x [1, 2)~.$
This is consistent with $r_0 \ua 1$ and thus the completions are given
by ${_{\RBbb}}M \cong {_{\RBbb}}M_2 \cong M_2$ and
${_{{\SBbb}^1}}M \cong {_{{\SBbb}^1}}M_2 \cong \Ga \bsl M_2$
with orbit space $[0, 2)~.$

Note that the $\RBbb$-- and the $\SBbb^1$--completion above are Hausdorff 
and essentially consist in completing the orbits in the gap created by the 
removed wedge in the annulus and then dividing by the discrete central 
subgroup $\Ga = \ZBbb_n$~. 
This is due to the fact that the imcomplete orbits in this example are still 
connected. 
If we remove in addition the half--open interval 
$1 \leq r < 2~,~\theta = \pi$~, the $\RBbb$-- and the $\SBbb^1$--completion 
of the resulting manifold will not be Hausdorff anymore, but the 
$\SBbb^1$--completion will retain the orbifold structure at the center.  

\section{Proper ${\g}\label{nmb:4}$--manifolds}

Let $(M, \zeta)$ be an orbifold--like $\g$--manifold.

\subsection{Definition}\label{nmb:4.1}~
The orbifold--like $\fr g$--manifold $(M, \zeta)$ is {\it proper} if
the $\wti G$--action on the generated $\wti{G}$--manifold
${_{\wti G}}M$ is proper, that is the mapping
\begin{equation*}
\bar \mu ~:~ {\wti G} \x {_{\wti G}}M \lra {_{\wti G}}M \x {_{\wti G}}M~,
\end{equation*}
given by $\mu (g, z) = (g z, z)$ is proper in the sense that compact sets
have compact inverse images (no assumption of Hausdorff).

\subsection{Sequential Definition}\label{nmb:4.2}

The ${\wti G}$--action on ${_{\wti G}}M$
is proper, iff the following holds. Given sequences $\{ g_n \}~, ~\{ z_n \}$
in  ${\wti G}$~, respectively ${_{\wti G}}M$~, such that the sequences
$\{ z_n \}~, ~\{ g_n z_n \}$ converge to $z, ~{\bar z}$ in ${_{\wti G}}M$~,
the sequence $\{ g_n \}$ in ${\wti G}$ has a convergent subsequence. 
%(accumulation point)~.
We will see below that the sequential definition involves convergent sequences 
only inside the charts $j_g : M \to V_g \subset {_{\wti G}}M$~, which are
Hausdorff as a consequence of our assumption.

\begin{prop}\label{nmb:4.3}~
The orbifold--like $\fr g$--manifold $(M, \ze)$ is proper, if and only if
the following condition is satisfied~: 

For convergent sequences $\{ x_n \} \to x~, ~\{ y_n \} \to y$ in $M~,$
a sequence $\{ g_n \in \wti{G} \}$ has a convergent subsequence 
$g \in {\wti G}$~,
provided that the following equivalent properties hold $~:$

\begin{enumerate}

\item[(1)]
$(e, y_n) \in L_{\zeta} (g_n, x_n) = g_n L_{\zeta} (e, x_n)$~, or
$(g_n, x_n) \in L_{\zeta} (e, y_n)$~,

%\noindent
%(this is P.M.'s definition)
\item[(2)]
$y_n \in I (e; g_n, x_n)$~, or
$x_n \in I (g_n; e, y_n) = I (e; g_n^{-1}, y_n)$~.

\end{enumerate}
\end{prop}

\begin{demo}{Proof}
Without loss of generality, we may assume that
$z = [e, x] \in V_e \subseteq {_{\wti G}}M$~.
Then $z_n \in V_e$~, that is the leaves represented by $z_n$ intersect the
transversal $i_e (M)$~, at least for sufficiently large $n$~, and we may
write $z_n = [e, x_n]$~, such that $\{ x_n \} \to x$ is convergent in $M$~.
The same applies to the sequence $\{ g_n z_n \} \to {\bar z}$~.
We write $\bar z = [g_0, y] \in V_{g_0}~, ~[e, y] \in V \subseteq {_{\wti G}}M$~.
Then $g_n z_n \in V_{g_0}$~, that is the leaves represented by $g_n z_n$
intersect the transversal $i_{g_0} (M)$~, at least for sufficiently large $n$~,
and we may write $g_n z_n = [g_0, y_n] = g_0 [e, y_n]$~, such that
$\{ y_n \} \to y$ is convergent in $M$~.
Finally, $g_n z_n$ may be written as
$g_n z_n = [g_0, y_n] = g_n ~[e, x_n] = [g_n, x_n]$~,
that is $(g_0, y_n) \in L_{\zeta} (g_n, x_n)$~, or equivalently
$(g_n, x_n) \in L_{\zeta} (g_0, y_n)$~.
In terms of intersection sets, this may also be formulated as
$y_n \in I (g_0; g_n, x_n) = I (e; g_0^{-1} g_n, x_n)$~, or equivalently
$x_n \in I (g_n; g_0, y_n) = I (e; g_n^{-1} g_0, y_n)$~.
Therefore the properness of the $\fr g$--manifold $M$
may be reformulated in terms of sequences as follows:

\noindent
If $\{ x_n \} \to x~, ~\{ y_n \} \to y$ are convergent sequences in $M$
and $\{ g_n \}$ is a sequence in ${\wti G}$~, then $\{ g_n \}$ has a convergent
subsequence $g \in {\wti G}$~, provided that the following
equivalent properties hold~:~

\begin{enumerate}

\item[(3)]
There exists $g_0 \in {\wti G}$~, such that
$(g_0, y_n) \in L_{\zeta} (g_n, x_n)$~, or $(g_n, x_n) \in L_{\zeta} (g_0, y_n)$~,

\item[(4)]
There exists $g_0 \in {\wti G}$~, such that
$y_n \in I (g_0; g_n, x_n) = I (e; g_0^{-1} g_n, x_n)$~, or
$x_n \in I (g_n; g_0, y_n) = I (e; g_n^{-1} g_0, y_n)$~.

\end{enumerate}

\noindent
For any convergent subsequence of $\{ g_n \} \to g$~, we have then
\begin{equation*}
\{ g_n z_n = [g_0, y_n] = [g_n, x_n]  \} \to
{\hat z} = [g_0, y] = [g, x] = g ~[e, x] = g z~.
\tag{5}
\end{equation*}
The necessity of the conditions in Proposition \ref{nmb:4.3} is now clear
by taking $g_0 = e~.$ Conversely, if the convergent sequence
$\{ g_n z_n \}$ is of the form $\{ g_n z_n = [g_0, y_n] \}~,$ the modified
sequence $\{ g_n' = g_0^{-1} g_n \}$ satisfies
$\{ g_n' z_n = g_0^{-1} ( g_n z_n ) = [e, y_n] \}$
and hence the conditions in Proposition \ref{nmb:4.3}. Therefore it has a
convergent  subsequence $\{ g_n' \} \to g'$~. Then the convergent
subsequence $\{ g_n = g_0 g_n' \} \to g = g_0 g'$ has the required
property.
\qed
\end{demo}

Since a proper $\fr g$--action on $M$ is discretely valued, the
Hausdorff separation property holds inside the open sets 
$V_g \subset {_{\wti G}}M$ and can fail only for pairs 
of limit elements by Proposition \ref{nmb:2.10}. 

\begin{lem}\label{nmb:4.4}~

\begin{enumerate}
\item[(1)]
If the $\fr g$--action on $M$ is proper, then ${_{\wti G}}M$ is a proper
$H$--manifold for any closed subgroup $H \subset {\wti G}$~.

\item[(2)]
If the $\fr g$--action on $M$ is proper, and $\Ga \subset {\wti G}$ is a
discrete subgroup, then $\Ga$ acts properly discontinuously on ${_{\wti G}}M$
and the isotropy groups $\Ga_z$ at $z \in {_{\wti G}}M$ are finite.

\item[(3)]
If the $\fr g$--action on $M$ is proper,
then the isotropy groups $G_z$ at $z = [g, x] \in {_{\wti G}}M$ are compact.

\end{enumerate}
\end{lem}

\begin{demo}{Proof}
\thetag{1} and \thetag{2} are immediate. 
\thetag{3} follows from a familiar argument. 
For any sequence $\{ g_n \}$ in $G_z$~, we have
$g_n [g, x] = [g_n g, x] = [g, x]$ and therefore $\{ g_n \}$ must have an
accumulation point in $G$~, which must be in $G_z$ since
$G_z \subset {\wti G}$ is closed. Equivalently, $G_z$ is given by
$G_z \cong G_z \x \{ z \} = \bar\mu^{-1} (z, z)$~, which is compact.
\qed
\end{demo}

\begin{prop}\label{nmb:4.5}~
Let $(M, \zeta)$ be a orbifold--like $\fr g$--manifold. 
Then the following holds: 

\begin{enumerate}

\item[(1)]
If the $\fr g$--action on $M$ is proper, then ${_{\wti G}}M$ is a 
smooth orbifold, the discrete
central subgroup $Z = \pi_1 (G) \subset \wti{G}$ acts properly
discontinuously on ${_{\wti G}}M$ and the orbit space
${_G}M = G \x_{\fr g} M \cong Z ~\bsl~ {_{\wti G}}M$
is a proper $G$--orbifold admitting
${_{\wti G}}M \to {_G}M = Z ~\bsl~ {_{\wti G}}M$
as a (ramified) Galois covering. 

\item[(2)]
Conversely, if the $G$--action on
${_G}M$ is proper, then the $\fr g$--action on $M~,$ that is the
${\wti G}$-action on ${_{\wti G}}M$ is proper, if and only if
$Z = \pi_1 (G) \subset \wti{G}$ acts properly discontinuously on
${_{\wti G}}M~.$

%\item[(3)]

\end{enumerate}
\end{prop}

\begin{demo}{Proof}
\thetag{1} follows from \thetag{\ref{nmb:3.3}.1}, since $Z \subset \wti{G}$ 
is a discrete central subgroup.  
\thetag{2} follows also from \thetag{\ref{nmb:3.3}.1} by a similar argument. 
\qed
\end{demo}

Thus the properness condition is strongest for the $\wti G$--orbifold 
${_{\wti G}}M$~.

\subsection{The slice theorem for proper ${\fr g}$--manifolds}\label{nmb:4.7}

We recall the slice theorem for proper $G$--actions
in a form convenient for our purposes cf. %~\Mi~.

For any $x \in M$~, there exists a submanifold (slice)
$x \in S_x \subset M$~, satisfying the following properties~:

\begin{enumerate}

\item[(1)]
$G_x \cdot S_x \subseteq S_x$~, that is $S_x$ is invariant
under the action of the stabilizer $G_x \subset G$ of $G$ at $x~;$

\item[(2)]
There exists a $G$--equivariant isomorphism
$G \x_{G_x} S_x \osetl{\cong} G \cdot S_x~,$ such that the diagram
\begin{equation*}
\xymatrix{G \x_{G_x} S_x    \ar[rr]^{{\cong}}_{}  \ar[dd]^{{\pi}}_{}  &  &  G \cdot S_x   \ar[dd]^{{}}_{} \\
 &  &         &  &   \\
G  / {G_x}        \ar[rr]^{{\cong}}_{}  &  &  G (x)}
\tag{3}
\end{equation*}
is commutative and the tube $G \cdot S_x \subset M$ is open in $M~;$

\item[(4)]
If $g \cdot S_x \cap S_x \neq \emptyset$ for $g \in G~,$ then $g \in G_x~.$

\end{enumerate}

\noindent
Note that \thetag{4} is a consequence of \thetag{1} and
\thetag{2}~.

\begin{thm}\label{nmb:4.8}~
Let $(M, \zeta^M)$ be a proper $\fr g$--manifold.
For any $x \in M$~, there exists a submanifold $(\text{slice})$
$S_x \subset M$~, satisfying the follwing properties$~:$

\begin{enumerate}

\item[(1)]
$\zeta ({\fr g}_x) \subset \eX (S_x)~;$

\item[(2)]
$\zeta ({\fr g}) (y) + T_y (S_x) = T_y (M)~, ~y \in S_x~;$ 

\item[(3)]
$\zeta ({\fr g}) (x) \oplus T_x (S_x) = T_x (M)~;$

\item[(4)]
The tube $\Ga (\fr g) (S_x) \subset M$ is open in $M~;$
%The tube $\zeta (\fr g) (S_x) \subset M$ is open in $M~;$

\item[(5)]
If $X \in {\fr g}$ and $\zeta_X (y) \in T_y (S_x)$ for some
$y \in S_x~,$ then $X \in {\fr g}_x~.$

\end{enumerate}
\end{thm}

\begin{demo}{Proof}
Since ${_{\wti G}}M = {\wti G} \x_{\fr g} M$ is a proper $\wti{G}$--manifold,
we may invoke the slice theorem for a proper $\wti{G}$--action and obtain
the above slice decomposition \thetag{\ref{nmb:4.7}.1} to \thetag{\ref{nmb:4.7}.4} for ${_{\wti G}}M~.$
This implies the infinitesimal properties \thetag{1} to
\thetag{5} for ${_{\wti G}}M~.$
For $x \in M$~, we may choose the $\wti{G}_{j (x)}$--invariant slice
$\wti{S}_{j (x)} \subset V$~, since the stabilizer subgroup $\wti{G}_{j (x)}$
is compact.
Since the chart $j : M \to V \subset {_{\wti G}}M$ is a $\fr g$--equivariant
local diffeomorphism onto $V$~, we obtain a ${\fr g}_x$--invariant slice
$S_x \subset M$ satisfying the infinitesimal properties
\thetag{1} to \thetag{5}~.
\qed
\end{demo}

\subsection{Existence of ${\fr g}$--invariant metrics}\label{nmb:4.9}

Given a Riemannian metric $h$ on $M$~, the following are equivalent~: ~

\begin{enumerate}

\item[(1)]
$\zeta : {\fr g} \to \eX (M)$ takes values in the
Lie algebra of Killing fields of $h$~, that is
$\zeta : {\fr g} \to \Ka_h (M) = \set{X \in \eX (M)}{L_X ~h = 0}$~;

\item[(2)]
The normal bundle $Q_{\folF_{\zeta}}$ of $\folF_{\zeta}$ admits a $G$--invariant
metric $h_Q$~, such that
\begin{equation*}
L_\xi ~h_Q = 0~, ~\forall \xi \in \Cinfty{T \folF_{\zeta}}~;
\end{equation*}

\item[(3)]
The holonomy pseudogroup of $\folF_{\zeta}$ acts by local isometries on $(M, h)$~.

\end{enumerate}

\noindent
The second statement says of course that the $G$--equivariant foliation
$(\folF_{\zeta}, h_Q)$ is a ($G$--equivariant) {\it Riemannian foliation}.

The Riemannian metric $h$ on $M$ is complete, if and only if
$(\folF_{\zeta}, h_Q)$ is transversally complete.

It is well known that the existence of a $G$--invariant Riemannian
metric implies the properness of a $G$-action on a manifold $M$~.
In our context this can be formulated as follows.

\begin{prop}\label{nmb:4.10}~
Suppose that the $\fr g$--action $\zeta$ on $M$ is effective and
that ${_G}M = G \x_{\fr g} M$ satisfies the Hausdorff property.
If there exists a $(\text{complete})$
Riemannian metric $h$ on $M$~, such that the injective Lie homomorphism
$\zeta : {\fr g} \to \eX (M)$ consists of Killing
fields relative to $h$~, then there exists a $(\text{complete})$
Riemannian metric  $\ti h$ on the leaf space ${_G}M = G \x_{\fr g} M$~,
such that the following conditions hold$~:~$

\begin{enumerate}

\item[(1)]
$G$ acts by isometries on the leaf space ${_G}M = G \x_{\fr g} M$~,
that is $G \subseteq \Iso_{\ti h} ({_G}M)~;$

\item[(2)]
The charts $j_g : M \to {_G}M$ are isometries of $\fr g$--manifolds.

\item[(3)]
The action of the closure $\wbar{G} \subseteq \Iso_{\ti h} ({_G}M)$ of
$G \subseteq \Iso_{\ti h} ({_G}M)$ is proper.

\end{enumerate}
\end{prop}

\begin{thm}\label{nmb:4.11}~
If the $\fr g$--action $\zeta$ on $M$ is proper, then there
exist $(\text{complete})$\ Riemannian metrics $h$ on $M$ and
$\ti h$ on ${_{\wti G}}M = {\wti G} \x_{\fr g} M $~, such that the
following conditions hold$~:~$

\begin{enumerate}

\item[(1)]
$\wti{G}$ acts by isometries on the leaf space
${_{\wti G}}M = {\wti G} \x_{\fr g} M$~,
that is ${\wti G} \to \Iso_{\ti h} ({_{\wti G}}M)~;$

\item[(2)]
The charts $j_g : M \to {_{\wti G}}M$ are isometries of $\fr g$--manifolds$;$

\item[(3)]
$\zeta : {\fr g} \to \Ka_h (M) \subset \eX (M)~.$

\end{enumerate}
\end{thm}

\begin{demo}{Proof}
The proof follows the by now familiar pattern by working on the proper
$\wti{G}$--manifold ${_{\wti G}}M = {\wti G} \x_{\fr g} M~.$ 
From the slice theorem for proper $\wti{G}$--actions %\Mi~,
we know that there exists a $\wti{G}$--invariant Riemannian metric
${\ti h}$ on ${_{\wti G}}M$~, so that \thetag{1} is satisfied.
Pulling back ${\ti h}$ to a metric on $h = j^* {\ti h}$ on $M$
via the local diffeomorphism $j : M \to V \subset {_{\wti G}}M~,$
property \thetag{2} is satisfied by $(1)$ and the transition
formula \thetag{\ref{nmb:2.3}.2}, that is
$j_g = L_g \circ j_e = L_g \circ j~.$
\thetag{3} is a direct consequence of \thetag{1} and
\thetag{2}. \qed
\end{demo}

%Can the proof be carried out by working on the $\fr g$--manifold
%$(M, \zeta^M)$ and using the Slice Theorem \ref{nmb:4.7} directly?

\begin{cor}\label{nmb:4.12}~
If ${_G}M \cong Z ~\bsl~ {_{\wti G}}M$ is a manifold and the $G$--action
on ${_G}M$ is proper, then ${_{\wti G}}M$ is a proper $\wti{G}$--manifold,
that is the $\fr g$--action on $M$ is proper.

\end{cor}

\subsection{The pseudogroup of local isometries}\label{nmb:4.13}

\section{$\fr g$--vector bundles and equivariant vector bundles}\label{nmb:5}

\subsection{}\label{nmb:5.1}
For a $G$--equivariant vector bundle $E \osetl{\pi} M$~, the associated 
infinitesimal action $X \bra \wti{X}^*, ~X \in {\g}$ is projectable to 
the corresponding vector field $X^*$ determined by the $G$--action on $M$~. 
%${\ell_g}_* \wti{X}^*_u = \wti{X}^*_{\ell_g (u)}$
Further, since the action $\ti{\ell}_g : E_u \to E_{\ell_g (u)}$ is linear, 
the vertical component of the vector field $\wti{X}^*$with respect to 
a local triviaalization of $E$ is a linear vector field on $E$~. This notion 
is clearly independent of the choice of local tricalization. 
Thus we define the Lie algebra $\eX_{\proj} (E)$ as the Lie algebra of 
vector fields $Y$ on $E$ satisfying the following conditions:

\begin{enumerate}

\item[(1)]
$Y$ is $\pi$--projectable;

\item[(2)]
The vertical component of $Y$ is a linear vector field.

\end{enumerate}

Let now $(M, \zeta)$ be a $\g$--manifold. A $\g$--structure on 
a vector bundle $E \osetl{\pi} M$ is given by a Lie algebra 
homomorphism $\ti\ze : \g \osetl{\pi} \eX_{\proj} (E)$~, such that 
$\ti\ze$ lifts the $\g$--structure on $M$~, that is the following   
diagram is commutative
\begin{equation*}
\xymatrix{
 & \eX_{\proj} (E)  \ar[dd]^{\pi_*}   \\
{\fr g} \ar@{-->}[ur]^{\ti \ze} \ar[dr]^{\ze}  \\
 & \eX (M).  }
\tag{3}
\end{equation*}
It is clear that the previous constructions, in partiular \ref{nmb:2.1}, 
\ref{nmb:2.3} apply to this situation as well. In particular, we obtain 
a $G$--equivariant foliation $\wti\folF_{\ti\ze}$ on $G \x E$ which 
projects under $\ti\pi = \id \x \pi$ to $G \x M$~, that is, $E$ carries 
the structure of a foliated vector bundle in the sense of \cite{5a}. 

\begin{prop}\label{nmb:5.2}~
Let $(\ti\ze, E)$ be a $\g$--structure on 
the vector bundle $E \osetl{\pi} M$~, satisfying condition  
\thetag{\thetag{\ref{nmb:5.1}.3}} above. Then 

\begin{enumerate}

\item[(1)]
There is a $G$--equivariant foliation $\wti\folF_{\ti\ze})$ on $G \x E$ 
which projects under $\ti\pi = \id \x \pi$ to $G \x M$~. Thisd defines 
a $G$--equivariant foliated bundle
\begin{equation*}
(G \x E, \wti\folF_{\ti\ze}) \osetl{\ti\pi} (G \x M, \folF_\ze)~. 
\end{equation*}

\item[(2)]
The $G$--completions 
\begin{equation*}
{_G}E \osetl{{_G}\pi} {_G}M
\end{equation*}
define a $G$--equivariant vector bundle. 

\end{enumerate}
\end{prop}

%\begin{demo}{Proof}
%\qed
%\end{demo}

\nocite{*}\bibliographystyle{plain}

\end{document}